\documentclass[10pt]{amsart}
\usepackage{amsmath,amsthm,amssymb,amscd}
\usepackage{epsfig}
\usepackage{eucal}

\begin{document}
\newcommand{\con} {{\rm config}}
\newcommand{\Con} {{\rm Config}}
\newcommand{\Ba} [1] {\ensuremath{{\mathcal {B}}^{#1}}}    
\newcommand{\Hu} {{\mathcal H}}                            
\newcommand{\ba}{{\mathfrak b}}                            
\newcommand{\sba}{{\mathfrak {sb}}}                        
\newcommand{\pba}{{\mathfrak {pb}}}                        
\newcommand{\D} {{\mathfrak D}}                            
\newcommand{\B} {{\bf B}}                                  
\newcommand{\M}  [1] {\ensuremath{{\overline{\mathcal M}}{^{#1}_0(\R)}}}   
\newcommand{\oM} [1] {\ensuremath{{\mathcal M}_{0}^{#1}(\R)}}              
\newcommand{\dM} [1] {\ensuremath{{\widetilde{\mathcal M}}_{0}^{#1}(\R)}}  
\newcommand{\cM} [1] {\ensuremath{{\mathcal M}_{0}^{#1}(\C)}}              
\newcommand{\CM} [1] {\ensuremath{{\overline{\mathcal M}}{^{#1}_0(\C)}}}   
\newcommand{\Sh} {{\mathcal{S}}}                           
\newcommand{\Pb} {{\bf P}}                                 
\newcommand{\C} {{\mathbb C}}                              
\newcommand{\R} {{\mathbb R}}                              
\newcommand{\Z} {{\mathbb Z}}                              
\newcommand{\N} {{\mathbb N}}                              
\newcommand{\Pj} {{\mathbb P}}                             
\newcommand{\T} {{\mathbb T}}                              
\newcommand{\Sg} {\mathbb S}                               
\newcommand{\sg} {\sigma}                                  
\newcommand{\Gl} {{\rm Gl}}                                
\newcommand{\G} {{\mathcal G}}                             
\newcommand{\IG} {{\mathfrak S}}                           
\newcommand{\Cox} [1] {\ensuremath{J_{#1}}}                
\newcommand{\qCox} [1] {\ensuremath{\tilde{J}_{#1}}}       
\newcommand{\Bnd} [2] {{\ensuremath{\Lambda (#1,#2)}}}
\newcommand{\Cobnd} [2] {{\ensuremath{\Lambda^n(#1,#2)}}}
\newcommand{\trh} {\kappa}                                 
\newcommand{\Op} [1] {{\mathcal {O}}(#1)}                  
\newcommand{\Lc} [1] {{\mathcal {C}}(#1)}                  
\newcommand{\Li} [1] {{\mathcal {I}}(#1)}                  
\newcommand{\SI} {\ensuremath{SI}}                         
\newcommand{\la}{{\langle}}                                
\newcommand{\ra}{{\rangle}}                                

\theoremstyle{plain}
\newtheorem{thm}{Theorem}[subsection]
\newtheorem{prop}[thm]{Proposition}
\newtheorem{cor}[thm]{Corollary}
\newtheorem{lem}[thm]{Lemma}
\newtheorem{conj}[thm]{Conjecture}

\theoremstyle{definition}
\newtheorem{defn}[thm]{Definition}
\newtheorem{exmp}[thm]{Example}

\theoremstyle{remark}
\newtheorem*{rem}{Remark}
\newtheorem*{hnote}{Historical Note}
\newtheorem*{nota}{Notation}
\newtheorem*{ack}{Acknowledgments}

\numberwithin{equation}{section}

\title {Cellular Structures determined by Polygons and Trees}

\author{Satyan L. Devadoss}
\author{Ronald C. Read}

\address{Department of Mathematics, Ohio State University, Columbus, Ohio}
\email{devadoss@math.ohio-state.edu}
\address{Department of Combinatorics and Optimization, University of Waterloo, Canada}
\email{rcread@math.uwaterloo.ca}

\begin{abstract}
The polytope structure of the associahedron is decomposed into two categories, types and classes.  The classification of types is related to integer partitions, whereas the classes present a new combinatorial problem.  We solve this, generalizing the work of \cite{r1}, and incorporate the results into properties of the moduli space \M{n} studied in \cite{dev}.  Connections are discussed with relation to classic combinatorial problems as well as to other sciences.
\end{abstract}

\maketitle

\baselineskip=15pt

\vspace{-.1in}

%
%

\section {The Real Moduli Space}

\subsection{}
The Riemann moduli space ${\mathcal M}_g^n$ of surfaces of genus $g$ with $n$ marked points has become a central object in mathematical physics.  Introduced in Algebraic Geometry, there is a natural compactification ${\overline{\mathcal M}}{_g^n}$ of these spaces; their importance was emphasized by Grothendieck in his famous \emph{Esquisse d'un programme} \cite{gro}. The special case when $g=0$, the space \CM{n} of $n$ punctures on the sphere $\C\Pj^1$, is a building block leading to higher genera. They play a crucial role in the theory of Gromov-Witten invariants and symplectic geometry. Furthermore, they appear in the work of Kontsevich on quantum cohomology, and are closely related to the operads of homotopy theory.  We look at the {\em real} points \M{n} of this space; these are the set of fixed points under complex conjugation.  This moduli space will provide the motivation to study certain combinatorial objects.

\begin{defn}
The real {\em Deligne-Mumford-Knudsen} moduli space \M{n} is a compactification of the configurations of $n$ labeled points on $\R \Pj^1$ quotiented by the action of $\Pj \Gl_2(\R)$.
\end{defn}

\begin{figure}[h]
\centering {\includegraphics {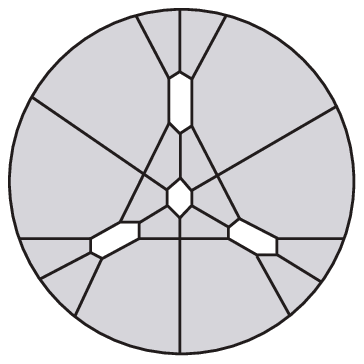}}
\caption{\M{5}}
\label{m05}
\end{figure}

\M{n} is a manifold without boundary of real dimension $n-3$. The moduli space is a point when $n=3$ and $\R\Pj^1$ when $n=4$. For $n>4$, these spaces become non-orientable. Figure~\ref{m05} shows \M{5} (shaded) as the connected sum of five copies of $\R \Pj^2$, with the four hexagonal `boundaries' and the outer circle carrying antipodal action.  In general, \M{n} can be described as the iterated blow-ups of real projective spaces (see \cite[\S4]{kap}, \cite[\S4]{dev}), getting extremely complicated as $n$ increases.

\begin{rem}
The diagram shown in Figure~\ref{m05} is first found in the work of Brahana and Coble in $1926$~\cite[\S1]{bc}, relating to possibilities of maps with twelve five-sided countries.
\end{rem}


\subsection{}
A beautiful fact about the moduli space is its tiling by the convex polytope known as the associahedron $K_n$.  It is the geometric realization of the poset of all meaningful bracketings on $n$ variables.  Originally defined for use in homotopy theory by Stasheff \cite[\S2]{jds}, it was later given a convex polytope realization by Lee \cite{lee}. We use an alternative definition:

\begin{defn}
The {\em associahedron} $K_n$ is a convex polytope of dim $n-2$ with codimension $k$ faces corresponding to using $k$ sets of non-intersecting diagonals on an $(n+1)$-gon.
\end{defn}

\begin{figure}[h]
\centering {\includegraphics {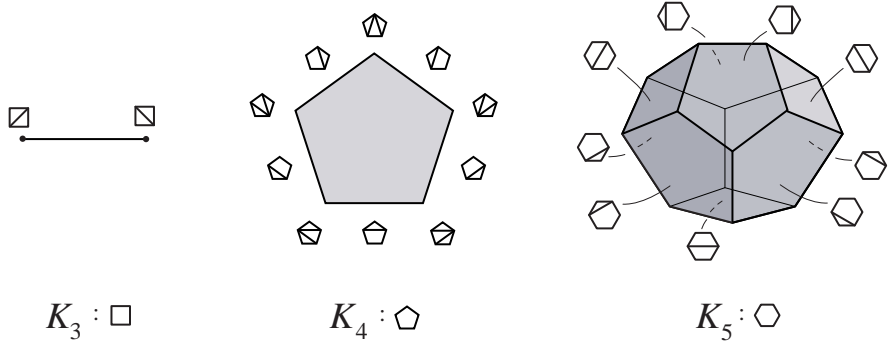}}
\caption{Associahedra $K_3$, $K_4$ and $K_5$}
\label{k3k4k5}
\end{figure}

\begin{exmp}
The associahedron $K_2$ is a point. Figure~\ref{k3k4k5} shows the associahedra $K_3$, $K_4$, and $K_5$. For instance, the polytope $K_4$ is a pentagon whose sides (codim one faces) are labeled with 5-gons with one diagonal and whose vertices (codim two faces) are labeled with 5-gons with two diagonals.\footnote{Mention of {\em diagonals} will henceforth mean non-intersecting ones.} 
\end{exmp}

\begin{figure}[h]
\centering {\includegraphics {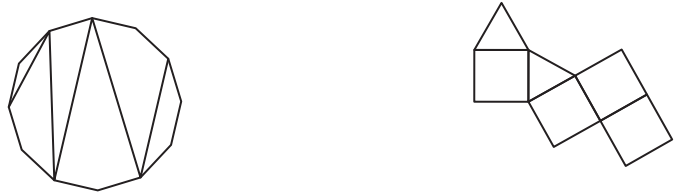}}
\caption{Polygons and clusters}
\label{cluster}
\end{figure}

We redraw our pictures of dissected polygons so that every region of the dissection becomes a regular polygon {\em without} diagonals. In the nomenclature of \cite{r1}, we denote the resulting object as a {\em cluster}, where each regular polygon of the cluster is called a {\em cell}. Thus the polygon shown in Figure~\ref{cluster}a with five diagonals is redrawn as a cluster in Figure~\ref{cluster}b having six cells. Clearly dissected polygons and clusters are equinumerous, and it is sometimes easier to work with clusters rather than with dissections, as the following important property of the associahedron demonstrates.

\begin{prop} \textup{\cite[\S2]{jds}} \label{p:product}
Each face of $K_n$ is a product of lower dimensional associahedra.
\end{prop}

\begin{figure} [h]
\centering {\includegraphics {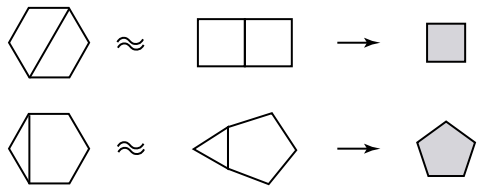}}
\caption{Codim one faces of $K_5$}
\label{k5codim1}
\end{figure}

\begin{exmp}
We look at the top (codim one) faces of $K_5$. The three-dimensional $K_5$ corresponds to a 6-gon, which has two distinct ways of adding a diagonal.  One way, in Figure~\ref{k5codim1}a, shows the $6$-gon as a cluster of two $4$-gon cells; each cell corresponds to the line segment $K_3$, where the product of two such is a square. Alternatively, Figure~\ref{k5codim1}b illustrates the other type of face, which is the product of $K_2$ (point) and $K_4$ (pentagon).
\end{exmp}

\subsection{}
It was shown in \cite[\S2]{dev} that \M{n} is tessellated by exactly $\frac{1}{2}(n-1)!$ copies of $K_{n-1}$. For example, observe from Figure~\ref{m05} that \M{5} is tiled by twelve $K_4$ pentagons.  The reasoning is as follows: There are $n!$ possible ways of labeling the sides of an $n$-gon, each corresponding to an associahedron.  However, since there exists a copy of the dihedral group $D_n$ of order $2n$ in $\Pj \Gl_2(\R)$, and since \M{n} is defined as a quotient by $\Pj \Gl_2(\R)$, two labeled polygons (associahedra) are identified by an action of $D_n$. Therefore, in some sense, each associahedral domain of \M{n} corresponds to a labeled $n$-gon up to rotation and reflection.

\begin{figure} [h]
\centering {\includegraphics {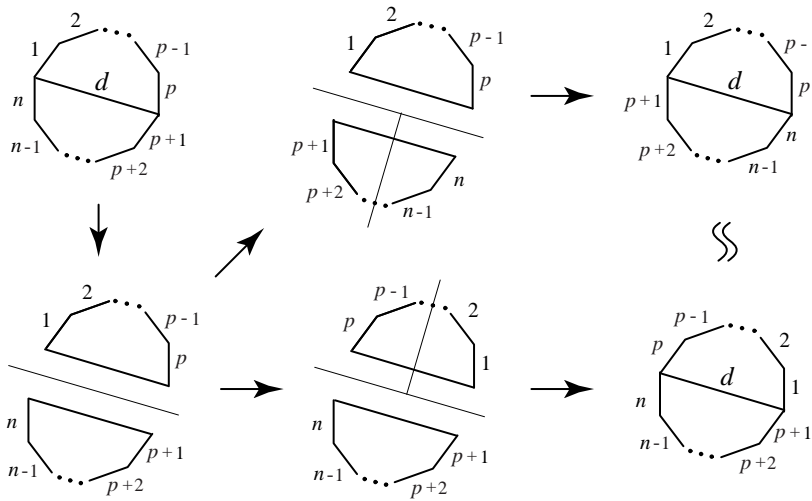}}
\caption{{\em Twist} along $d$}
\label{twist}
\end{figure}

\begin{defn}
Let $\G(n,k)$ be the set of equivalence classes of {\em labeled} regular $n$-gons with $k$ diagonals up to the action of $D_n$. Let $G$ be an element of $\G(n,k)$ and $d$ be a diagonal of $G$.  The {\em twist} of $G$ along $d$, denoted by $\nabla_d(G)$, is the polygon in $\G(n,k)$ obtained by separating $G$ along $d$ into two parts, `twisting' (reflecting) one of the pieces, and gluing them back (Figure~\ref{twist}).  It does not matter which piece is twisted since the two results are identified by an action of $D_n$.
\end{defn}

\begin{thm} \textup{\cite[\S4]{dev}}
Two polygons, $G_1, G_2 \in \G(n,k)$, representing codim $k$ faces of associahedra, are identified in $\M{n}$ if there exist diagonals $d_1, \ldots, d_r$ of $G_1$ such that $$(\nabla_{d_1} \cdots \nabla_{d_r}) (G_1) = G_2.$$
\label{t:glue}
\end{thm}

\begin{figure} [h]
\centering {\includegraphics {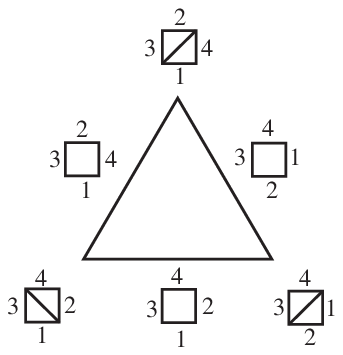}}
\caption{\M{4}}
\label{m04}
\end{figure}

\noindent The previous statement shows how copies of associahedron glue together
to form the moduli space. Figure~\ref{m04} demonstrates \M{4} tiled by three labeled $K_3$'s.

%
%

\section {Categories of Faces}

\subsection{}
Recall that faces of associahedra can be represented by dissected polygons.  We give the following definitions that categorize the faces of $K_n$.

\begin{defn}
Two polygons $G_1$ and $G_2$ are of the same \ldots 

\emph{dimension} if their corresponding faces are of the same dimension.

\emph{type} if their corresponding faces are the same polytope.

\emph{class} if they are identified under the actions of $D_n$ and twisting.

\noindent Dimension, type and class are in increasing order of finer partitions.
\end{defn}

\begin{figure} [h]
\centering {\includegraphics {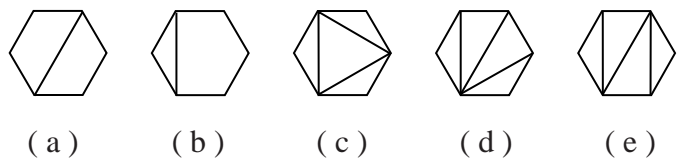}}
\caption{Different faces of $K_5$}
\label{types}
\end{figure}

\begin{exmp}
In Figure~\ref{types}, we note the difference between the categories of polygons belonging to $K_5$. Recall that the number of diagonals of a polygon is equal to the codim of its face; polygons (a) and (b) are of dim two, and (c) through (e) are of dim zero.  In particular, polygons (c) through (e) are of the same type (points) whereas (a) and (b) are not; the polytope (a) is a square and (b) a pentagon. Only (d) and (e) are of the same class; that is, twisting along a diagonal of (d) will result in (e). 
\end{exmp}

It makes sense to decompose the faces into dimensions and polytopes, but what does it mean to talk about classes?  The two different classes of points discussed above correspond to the following:  At one class (c), three pentagons meet; there are two such points in $K_5$.  The other class (d) is where two pentagons and one square meet; there are twelve such points.  Therefore, classes categorize {\em the different ways polytopes intersect} in $K_n$. This idea can be extended to better understand the structure of \M{n}.


\subsection{}
The combinatorics involving dimensions are handled in \cite{dev}.  We recall

\begin{lem}
The number of faces in $K_n$ of codim $k$ is $$\frac{1}{k+1} \; \binom{n+k}{k} \; \binom{n-2}{k}.$$ 
In other words, given the $f$-vector of the polytope $K_n$ to be $f=(f_0, \ldots, f_{n-3})$, the number above is $f_{n-k-2}$.
\label{l:cayley}
\end{lem}

\begin{proof}
This was originally proved by Professor A.\ Cayley in 1891 in finding the number of partitions of an $(n+1)$-gon into $k$ parts  \cite{cay}. The recent work of Przytycki and Sikora gives a short and elementary proof while providing relations to knot theory \cite{ps}. 
\end{proof}

\begin{rem}
The number of faces in $K_n$ of codim $n-2$ is the well-known Catalan number $c_n = \frac{1}{n+1} \binom{2n}{2}$.  That is, there are $c_n$ vertices in the associahedron $K_n$.  It is interesting to note that Stanley shows sixty-six ways to represent the Catalan numbers \cite[Exercise 6.19]{s2}.
\qed
\end{rem}

\noindent These are the three fundamental problems addressed in this paper:
\vspace{.1in}

{\em P1.} \   How many different types/classes  of codim $k$ are there in $K_n$ ?
\vspace{.03in}

{\em P2.} \   What are the different types/classes in $K_n$ of codim $k$? 
\vspace{.03in}

{\em P3.} \  How many faces of $K_n$ belong to a given type/class?
\vspace{.1in}

In some sense, faces of the associahedron correspond to polygons and faces of the moduli space to {\em labeled} polygons.  Therefore, since types and classes were defined regardless of labeling, the answers to {\em P1}\, and {\em P2}\, both carry over to the moduli space. This is not the case for {\em P3}\,. However, we find in \cite{dev} a formula (that follows from Theorem~\ref{t:glue}) which addresses this dilemma.

\begin{prop} 
Let $\Lambda(k, {\mathfrak X})$ denote the number of codim $k$ faces in the CW-complex ${\mathfrak X}$. Then
$$\ \Lambda(k, \M{n+1})  \ = \ {\displaystyle \frac{n!}{2^{k+1}}} \, \cdot \, \Lambda(k, K_n).$$
\label{p:k-to-m}
\end{prop}


\subsection{}
In what follows, we find solutions to our fundamental problems regarding the classification of types. We try to clarify the situation with the following example.

\begin{figure} [h]
\centering {\includegraphics {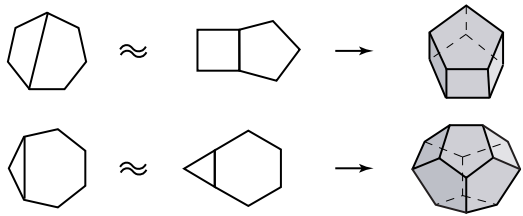}}
\caption{Codim one faces of $K_6$}
\label{k6codim1}
\end{figure}

\begin{exmp}
Looking at the case of $K_6$, Lemma~\ref{l:cayley} shows there to be $14$ codim one faces.  Figure~\ref{k6codim1} shows there to be two different types of such faces (answers {\em P1}\,), one being $K_6$ and the other a pentagonal prism (answers {\em P2}\,).  Out of the 14 faces, there turn out to be seven faces of each type (answers {\em P3}\,). 
\end{exmp}

Polygons are of the same type if their corresponding associahedral faces have identical polytopal structure. As Proposition~\ref{p:product} indicates, the polytope structure is solely determined by the number of different cells, not the way these cells are glued to form the cluster. The following notation is well-defined:  A face is of type $\la 3^{m_3} : 4^{m_4} : \cdots \ra$ if it corresponds to a cluster that has $m_i$~number of $i$-sided cells.  The superscript $m_i = 1$ as well as the terms with $m_i = 0$ are omitted. For example, the three types in Figure~\ref{types} are, respectively, $$\la 3^0 : 4^2 \ra = \la 4^2 \ra \hspace{.5in} \la 3^1 : 4^0 : 5^1 \ra = \la 3 : 5 \ra \hspace{.5in} \la 3^4 \ra.$$

Define $p_k(n)$ to be the number of partitions of $n \in \N$ into exactly $k$ integers $n_1, \ldots, n_k$, where $n_i > 0$ and $n = n_1 + \cdots + n_k$.  There is no closed form for $p_k(n)$ but the following recurrence relation provides a means to calculate it:
\begin{equation}
p_k(n) = p_{k-1}(n-1) + p_k(n-k).
\label{e:part}
\end{equation}

\begin{thm} [Solution to {\em P1}\, and {\em P2}\,]
The number of different types in $K_n$ of codim $k$ is $p_{k+1}(n-1)$. The different partitions correspond to the different types, where $n-1 = n_1 + \cdots + n_{k+1}$ is associated with the polytope $K_{n_1} \times \cdots \times K_{n_{k+1}}$.
\end{thm}

\begin{proof}
An $(n+1)$-gon with $k$ diagonals (that is, a cluster containing $k+1$ cells) represents a codim $k$ face of $K_n$. The type of a face depends on the polytope structure; this is soley determined by the cells forming the cluster, {\em not} the way in which they are arranged to form the $(n+1)$-gon. Associating an $(n+1)$-gon to the integer $n-1$, observe that partitioning the polygon into $k+1$ cells is identical to partitioning the integer into $k+1$ parts.
\end{proof}

\begin{thm} [Solution to {\em P3}\,]
When $k = -1 + \sum m_i$ and $n = 1 + \sum (i-2)m_i$, the number of faces in $K_n$ of type $\la 3^{m_3} : 4^{m_4} : \cdots \ra$ is 
$$\frac{1}{k+1} \; \binom{n+k}{k} \; \frac{(m_3 + m_4 + \cdots)!}{m_3! \ m_4! \ \cdots}.$$
\end{thm}

\begin{proof}
This is proven by Goulden and Jackson \cite[Exercise 2.7.14]{gj}. Their argument is similar to the proof given below for the enumeration of $A$-clusters (see~\S\ref{ss:aclusters}), but without the twist operation and with the extra variables $m_i$ in the generating function to mark the numbers of $i$-sided regions. From the recursive equation for the generating function, the coefficients (functions of $m_i$) are obtained by Lagrange's theorem and then expanded using the multinomial theorem.
\end{proof}

Values computed are shown in Table~\ref{t:diss} using Equation~\ref{e:part}.  For each size of polygon, and each distribution of sizes of the cells in the dissection, it gives the number of polygons which have that distribution.  An example of an entry for the partition of a 10-gon shows that there are 660 polygons of type $\la 3^2:4:6 \ra$.  That is, there are 660 faces in $K_9$ that have the polyhedral structure of $K_2 \times K_2 \times K_3 \times K_5$.

%
%

\section {Classes and the Combinatorics of Twisting}

\subsection{} \label{ss:aclusters}
Unlike types, the problem with regard to classes is no longer classical. Counting the number of dissected polygons up to rotation and reflection was studied in \cite{r1}.\footnote{This is equinumerous to the set of {\em outerplanar graphs}.}   The difference now comes from imposing the additional equivalence relation of twisting.  The method used to enumerate classes naturally follows as an extension of \cite{r1}.   

\begin{figure}[h]
\centering {\includegraphics {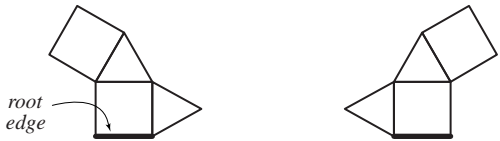}}
\caption{Outside root edge}
\label{rootedge}
\end{figure}

\begin{defn}
A cluster is said to be {\em rooted at an outside edge} if one edge (the {\em root}) of its perimeter is distinguished from the others.  Figure~\ref{rootedge}a shows the root edge as a thicker line. For convenience, we refer to clusters rooted at an outside edge as {\em $A$-clusters}.
\end{defn}

Let $a_{m,n}$ be the number of $A$-clusters having $m$ cells and $n$ outside edges not counting the root edge.   Define the generating function
\begin{equation}
A(x, y)  =  \sum a_{m,n}x^m y^n,
\label{e:gena}
\end{equation}
where for the moment the limits for $m$ and $n$ are unspecified. For an $A$-cluster, there will be one cell with the root edge on its boundary;  we can call this the {\em root cell}.  Assume this cell has $k$ unrooted edges.  At each of these edges we have the option of gluing the root of some $A$-cluster to it.\footnote{This gluing is identical to the composition map for operads in the context of rooted trees.} Each of these glued clusters contributes its values of $m$ and $n$ to those of the resulting $A$-cluster.  For the edges of the root cell at which nothing was glued, the contributions to $m$ and $n$ are 0 and 1 respectively. Hence the generating function $A(x,y)$ will contain the term $y$.  In the rest of $A(x,y)$, $m$ goes from 1 upwards and $n$ from $2$ upwards.

We now have a typical P\'{o}lya-class problem.\footnote{For information on such problems, see \cite{pol}, \cite{pr}, \cite{r3}.} There are $k$ {\em sites} around the root cell at each of which we can place a figure, an $A$-cluster or just an edge, and the parameters $m$ and $n$ are additive.  The figure counting series is the function $A(x,y)$ just defined.   To accommodate the twisting operation, we must regard an $A$-cluster as being equivalent to the one obtained by reversing the root edge --- as in Figure~\ref{rootedge}a and \ref{rootedge}b. This reversal interchanges the sites in pairs if $k$ is even;  if $k$ is odd then one site remains fixed while the others interchange in pairs.   Hence the cycle-index for this problem is
\begin{eqnarray*}
\frac{1}{2}(s_1^{2\alpha} + s_2^\alpha) & \hspace{1in} & k = 2\alpha \\
\frac{1}{2}s_1(s_1^{2\alpha} + s_2^\alpha) & \hspace{1in} & k = 2\alpha + 1. \\
\end{eqnarray*}

\vspace{-.2in}
\noindent By P\'{o}lya's theorem, the generating function for $A$-clusters for which the root cell has $k$ unrooted edges is obtained by putting $s_1 = A(x, y) \mbox{ and } s_2 = A(x^2, y^2)$ in the appropriate formula above. Summing for $k$ from 2 upwards we get a function of $s_1$ and $s_2$,
$$\frac{1}{2}\sum_{k\geq 2} s_1^k + \frac{1}{2}(1 + s_1)\sum_{\alpha \geq 1} s_2^\alpha,$$
which reduces to
$$\frac{1}{2} \left[ \frac{s_1^2}{1 - s_1} + \frac{(1+s_1)s_2}{1 - s_2} \right].$$
By substituting $s_1 = A(x, y) \mbox{ and } s_2 = A(x^2, y^2)$ in this expression, multiplying by $x$ to account for the root cell, and adding the term $y$ which we know must be present, we recover the generaing function $A(x,y)$ for $A$-clusters. This gives the recursive equation
\begin{equation}
A(x,y) = y + \frac{1}{2} \left[ \frac{A(x,y)^2}{1-A(x,y)} + 
\frac{(1+A(x,y))A(x^2,y^2)}{1-A(x^2,y^2)} \right].
\label{e:axy}
\end{equation}
The occurrence of $A(x^2, y^2)$ in this equation makes it impossible to obtain a closed form for the coefficients $a_{m,n}$ but they can be successively computed as far as one wishes (see Table~\ref{a-mn}).

\begin{rem}
Another proof of Lemma~\ref{l:cayley} is obtained by computing the corresponding result above when twisting is {\em not} allowed.  Then the relevant permutation group is the identity group, with cycle-index $s_1^k$, and the analog of Equation~\ref{e:axy} is
$$V(x,y) = \sum v_{m,n} x^m y^n  =  y + \frac{xV^2(x,y)}{1 - V(x,y)},$$
where the numbers $v_{m,n}$ are the counterparts to the $a_{m,n}$ under the new conditions. In this case, as was shown in~\cite{r1}, a closed result is possible:
$$v_{m,n}  = \frac{1}{m} \binom{n-2}{m-1} \binom{m+n-1}{n}.$$
That is, $v_{m,n}$ is the number of codim $m-1$ faces in $K_n$; values are given in Table~\ref{v-mn}.
\end{rem}


\subsection{}
\begin{figure}[h]
\centering {\includegraphics {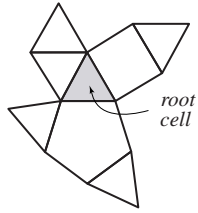}}
\caption{Root cell}
\label{rootcell}
\end{figure}

We now turn to clusters rooted at a cell, denoted as $B$-clusters. Figure~\ref{rootcell} shows a typical $B$-cluster, having a distinguished cell (shown shaded) called the {\em root cell}.  A $B$-cluster is obtained by gluing $A$-clusters (including possibly the empty $A$-cluster) by their roots to the edges of the root cell. If the root cell is $k$-sided, we have a P\'{o}lya-class problem for which the figure counting series is $A(x, y)$ and the permutation
group is the dihedral group $D_k$.

Since the twisting operation is allowed, the way that the $A$-clusters are attached is not important.  For this reason we are spared the tiresome necessity of distinguishing between symmetric and asymmetric clusters which had greatly complicated the enumerations performed in \cite{r1}. The cycle-index of $D_k$ is
\begin{eqnarray*}
\frac{1}{2k}\sum_{r|k}\varphi(r)s_r^{k/r}+ \frac{1}{4}(s_1^2s_2^{\alpha-1}+s_2^\alpha) & \hspace{.8in} & k = 2\alpha \\
\frac{1}{2k}\sum_{r|k}\varphi(r)s_r^{k/r}+\frac{1}{2}s_1s_2^\alpha & \hspace{.8in} & k = 2\alpha + 1. \\
\end{eqnarray*}

\vspace{-.2in}
\noindent The $B$-clusters for a given value of $k$ are enumerated by substituting $A(x, y)$ in this cycle-index to get the generating function denoted by $Z(D_k, A(x,y))$.   Hence the set of all $B$-clusters is enumerated by
$$B(x, y)  =  \sum b_{m,n} x^m y^n = x \sum_{k\geq 3} Z(D_k, A(x, y)).$$
There seems to be no way of simplifying this expression, but the coefficients can be easily computed as far as desired.  Some values for $b_{m,n}$ are given in Table~\ref{b-mn}.


\subsection{}
\begin{figure}[h]
\centering {\includegraphics {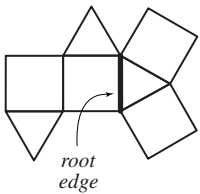}}
\caption{Inside root edge}
\label{irootedge}
\end{figure}

The next enumeration is for clusters rooted at an inside edge which are asymmetrical about that edge; that is, for which the two clusters at the edge are different (Figure~\ref{irootedge}).  Clearly these $C$-clusters, as they can be called, are formed by gluing together two different $A$-clusters, neither of which is the empty cluster, at their root edges. 

P\'{o}lya's theorem shows the generating function for these clusters would be obtained by substituting $A(x, y) - y$ in the cycle-index of the symmetric group $\Sg_2$, namely $\frac{1}{2}(s_1^2  +  s_2)$, but for the fact that we want the two $A$-clusters to be different.  A variant of P\'{o}lya's theorem, which is well-known but in any case easily verified, provides the required result by substituting instead in the expression $\frac{1}{2}(s_1^2  -  s_2)$. Calling the resulting generating function
$C(x, y)$, we have 
$$C(x, y) =\frac{1}{2}[(A(x, y)-y)^2 - (A(x^2, y^2)-y^2)].$$
The final step, going from the generating function for cell-rooted clusters to {\em free} clusters (those not rooted at all), is carried out by means of a theorem originally due to Otter \cite{ott} (also see \cite{hp}, \cite[\S4]{r1}). It states that the generating function for these free clusters, call it $F(x, y)$, is obtained by subtracting $C(x, y)$ from $B(x, y)$.  The number $f_{m,n}$ of clusters with $m$ cells and $n$ outside edges, is found from
$$F(x, y) = B(x,y) - C(x,y) = \sum f_{m,n}x^m y^n.$$
Some values for $f_{m,n}$ are given in Table~\ref{f-mn}.  We have just proved

\begin{thm} [Solution to {\em P1}\,]
The number of different classes of $K_n$ of codim $k$ is $f_{k+1,n+1}$. 
\end{thm}

%
%

\section {Trees and the Complex Moduli Space}

\subsection{}
\begin{figure}[h]
\centering {\includegraphics {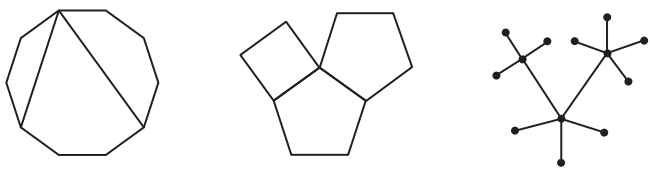}}
\caption{Polygons, clusters, and trees}
\label{f:pct}
\end{figure}

As a help in the study of the various classes, we prepared a set of drawings of the corresponding dissected polygons up to ten vertices. Figures 22-26 illustrate the different classes, labeled as $n.k.i$ to denote an $n$-gon with $k$ diagonals where $i$ ranges from one to $f_{n+1,k+1}$. The computer algorithm for doing this was nothing out of the ordinary, but a few comments on how the drawings were produced may be of some interest. We exploited the connection between dissected polygons and trees as illustrated in Figure~\ref{f:pct}.\footnote{Considering dissected polygons as plane graphs, we obtain trees by taking their {\em weak dual}.}
Trees dual to dissected polygons are {\em homeomorphically irreducible}, having no vertices of degree two.\footnote{They are also known as {\em topological} or {\em reduced} trees.} Henceforth, mention of trees will imply homeomorphically irreducible ones. To each dissected polygon, we associate a tree imbedded in the plane. For each diagonal of the polygon, there is a tree edge connecting the two vertices on either side. The twist operation corresponds to reflecting one of the two subtrees at the ends of a given edge. Classes correspond to planar imbeddings of trees up to twisting.  In other words, for a given class, each vertex of the tree is given the addition information of preserving the cyclic order of edges around that vertex. The problem was therefore to produce one such tree for each class.

Since the trees needed were quite small, it was not worth devising a clever algorithm; both elegance and efficiency were sacrificed for ease of programming.  Rooted trees, the root being a vertex of the tree, were generated step by step according to height.   There is one tree of height zero, and those of height one are the stars\footnote{Also known as a {\em corolla}, a star is a collection of edges meeting at a common vertex.} with three or more edges at the root.  Those of height two were obtained by grafting a subtree of height at most one, by its root, to the terminal vertices of these stars.  In general, trees of height $k$ were obtained by grafting trees of height at most $k - 1$. At each stage, if two trees were reflections of each other, only one was retained, thus taking care of the twist operation.  This corresponds to producing clusters rooted at an edge, which made possible a check that the right numbers of trees were being produced.

The resulting lists of rooted trees would usually contain several copies of the same unrooted tree.   The unwanted duplicates were weeded out by computing for each tree a {\em code} similar to that described in \cite{r2} and eliminating all but one of those with the same code.   The actual program contained some refinements to the procedure just described but we shall skip the details. It was a fairly simple task to go from trees to dissected polygons.  The diagrams\footnote{Although diagrams were created for ten-sided polygons, we do not include them here for the sake of space.} were then produced automatically using programs that had already been developed for the production of the book {\em An Atlas of Graphs} \cite{rw}.

\begin{prop}[Solution to {\em P2}\,]
The different classes correspond to the different pictures in Figures 22 - 26.
\end{prop}

At a later stage we decided to compute for each class the number of fixed dissections of that class.   Had we thought of this earlier, it could have been done more easily along with the generation of the classes, but as it was, we derived these numbers directly from the data on the dissected polygons.   For each such polygon, we produced the fixed polygons obtainable from it by rotations or by combinations of the twist operation.   Note that twisting about some chords leaves the diagram unchanged, in which case there is no point in performing the twist.   For example, those chords which have a single polygon on one side fall into this category, and are easily recognised.   For ease of programming, only the chords of this kind were exempted from the twist operation; no attempt was made to recognize any others. As a result, duplicates could, and were, produced.   However, for polygons with at most ten sides there can be at most five such chords, so the number of duplicates was small.   They were eliminated by a suitable coding procedure, similar to that mentioned above for trees. It is reassuring that by summing these numbers for the appropriate sets of classes we recover the numbers given in Table~\ref{t:diss}.

\begin{prop}[Solution to {\em P3}\,]
The number of faces of $K_n$ of a particular class is indicated in the top right-hand corner of each picture in Figures 22 - 26.
\label{p:cp3}
\end{prop}

\subsection{}  As mentioned above, there is a duality between clusters and trees.  Clearly each cell of a cluster corresonds to a star in the tree. A face of the associahedron of a given type only recognizes the set of stars, not the way in which they are grafted to form the tree. However, a given class not only sees the underlying tree, but the additional cyclic order of edges around each vertex. Hence, faces of associahedra can be categorized by types, trees, and classes, in increasing order of refinement.

\begin{figure}[h]
\centering {\includegraphics {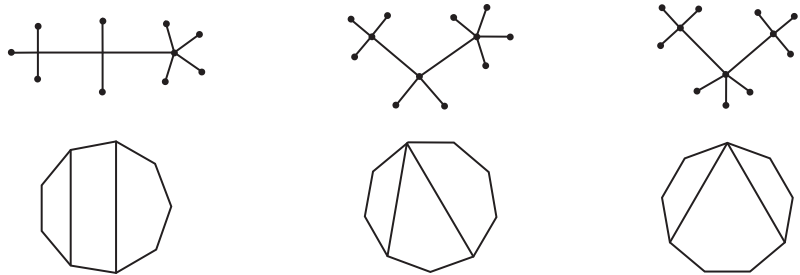}}
\caption{Some faces of $K_8$ and their trees}
\label{phylo}
\end{figure}

\begin{exmp}
Figure~\ref{phylo} shows three different trees and their corresponding (dual) polygons.  All of them are of the same type $[4^2 : 5]$ and all fall into different classes.  Note, however, that the first two correspond to the same tree, whereas the third does not.  By not preserving the cyclic order of edges, we obtain a distinction that falls between the categories of type and class.
\end{exmp}

The enumeration of all trees (not just homeomorphically irreducible ones) by their degree partitions is given by Harary and Prins \cite[\S4]{hpr}.  They provide a method for computing the generating function
$$\overline{H}(x, t_1, t_2, \ldots) = \sum \overline{h} \, x^m t_1^{m_1} t_2^{m_2} \ldots,$$
where $\overline{h}$ is the number of trees having $n$ vertices altogether, with $n_i$ vertices of degree $i$.  To get the number of {\em homeomorphically irreducible} trees $h$ partitioned by degree, we remove the possibility of any degree being two (letting $t_2=0$). This involves only standard enumerative procedures and the details are omitted.  We give below the equations that determine the required generating function for homeomorphically irreducible trees.\footnote{The symmetric function notation for cycle-indexes has been used (see \cite{r0}). For example, we write $h_n[P]$ in place of $Z(\Sg_n, P(x, t_1, t_2, \ldots))$.}
\begin{eqnarray*}
P &=& x(t_1 + t_3 h_2[P] + t_4 h_3[P] + \ldots ) \\
R &=& x(t_1 P + t_3 h_3[P] + t_4 h_4[P] + \ldots )\\
H &=& R - a_2[P]
\end{eqnarray*}

\begin{thm}
The number of different trees in $K_n$ of codim $k$ is given by the above generating function
$$H = H(x, t_1, t_3, \ldots) = \sum h \; x^m t_1^{m_1} t_3^{m_3} \ldots,$$ where $m=n+k+1$, $m_1=n$, and $\displaystyle{\sum_{i \geq 3} m_i = k+1}$.
\end{thm}

\noindent We note in passing that all homeomorphically irreducible trees up to $16$ vertices are depicted in \cite{rw} together with their degree partitions.


\subsection{}
\begin{figure}[h]
\centering {\includegraphics {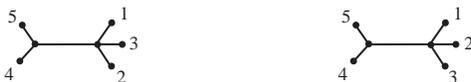}}
\caption{Identical labeled trees}
\label{phylolabel}
\end{figure}

Roughly speaking, faces of the associahedron correspond to polygons and faces of the moduli space to labeled polygons. Both types (polytope structures) and classes (intersections of polytopes) are natural objects in the associahedron and the moduli space.  We extend our idea of trees by labeling their external vertices.  Unfortunately, {\em labeled trees} are of little use to \M{n}; the trees shown in Figure~\ref{phylolabel} are identical considered as labeled trees and yet distinct as faces in \M{n}.
 
However, labeled trees are {\em exactly} what are needed to understand the points of the {\em complex} moduli space \CM{n}.  A rough sketch for this is as follows: The real moduli space \M{n} comes from a configuration of labeled points on $\R \Pj^1$ (viewed as a circle).  Keeping track of the cyclic ordering of these points is crucial, since two adjacent labeled points can switch places only by  first colliding with each other.  Points on a circle give rise to the planarity of trees, and the collision/permutation of adjacent points are kept track of with the twist operation. The situation is different for the complex case.  \CM{n} is based on a configuration of labeled points on $\C \Pj^1$ (viewed as a sphere).  There are diffeomorphisms of the sphere which can permute any two labeled points without their need to collide, giving rise to non-planar trees. This is the fundamental reason that the underlying structures are classes for \M{n} and trees for \CM{n}.

\begin{rem} We point out a few observations about moduli spaces and trees.

1. The codim $k$ faces of \CM{n} correspond to labeled trees with $n$ external vertices and $k+1$ internal vertices. In particular, both \M{n} and \CM{n} have the same number of zero dim cells, enumerated by binary trees with labeled external vertices.  In other words, the twist operation on planar binary trees removes the restriction of planarity.

2. There are operad structures related to the moduli spaces \M{n} \cite{dev} and \CM{n} \cite{gk}, classically  described in terms of labeled trees \cite[\S 1.4]{bv}. The associahedron is related to $A_\infty$ structures, whereas \M{n} seems to play a similar role in Fukaya's Lagrangian Floer cohomology.

3. Kapranov \cite{kap} has defined a natural double cover of the moduli space \M{n} by fixing one of the $n$ distinct labeled points in $\R \Pj^1$ to be $\infty$. Similar to \M{n}, this double cover also has an underlying structure of labeled trees \cite[\S 4.3]{dev} with two exceptions. First, the trees are no longer free but rooted at the label $\infty$.  Second, twisting is allowed except along the edge containing the root.

4. The reason for these moduli spaces to be characterized by trees without degree two vertices is due to the stability condition coming from Geometric Invariant Theory \cite[\S 8]{git}.  GIT gives a natural compactification of the space of  genus zero algebraic curves which are {\em stable} in the sense of having only finitely many automorphisms; that is, the curves cannot have just two distinguished points.
\qed

\end{rem}

%
%

\section {The Isotropy Group}

\subsection{}
There is a natural action of the symmetric group $\Sg_n$ on \M{n} which permutes the labelings of the $n$-gon. This action is certainly not free; therefore, for a given face $f$ in \M{n}, its isotropy group $\IG_f \subset \Sg_n$ need not be trivial, where  
$$\IG_f = \{\sg \in \Sg_n \,|\, \sg(f) = f\}.$$
Given two distinct faces of \M{n} which are of the same class, their isotropy groups are {\em isomorphic}; in fact, they are conjugate subgroups in $\Sg_n$.  We wish to determine, up to isomorphism, the isotropy group of a given class; we interchangably write $\IG_f$ to denote the isotropy group of $f$ and its isomorphism class.

\begin{figure}[h]
\centering {\includegraphics {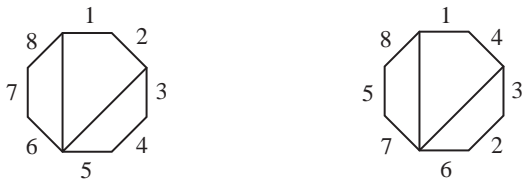}}
\caption{Faces with conjugate subgroups}
\label{isotropy}
\end{figure}

\begin{exmp}
Let $f$ be the face of \M{8} given in Figure~\ref{isotropy}a; its class is $8.2.7$.  Its isotropy group is generated by the reflections $\{(68),\, (35),\, (12)(83)(74)(65)\}$. Since the reflections commute,  $\IG_f$ is isomorphic to the direct product $\Z_2^3$. As discussed above, all isotropy groups of class $8.2.7$ will be isomorphic.  In particular, a face $g$ of the same class, as shown in Figure~\ref{isotropy}b, will have a conjugate subgroup  $\IG_g = (24)(567) \cdot \IG_f \cdot (765)(42)$.
\label{e:isotropy}
\end{exmp}

\begin{cor}
For a given class $n.k.i$, let $\trh$ be the number of faces in $K_{n-1}$ of this class.  The number of faces in \M{n} of the class $n.k.i$ is
$$\trh \cdot \frac{(n-1)!}{2^{k+1}}.$$
\end{cor}

\begin{proof}
This immediately follows from Propositions~\ref{p:k-to-m} and \ref{p:cp3}.
\end{proof}

This can be obtained using another approach:  Choose a dissected polygon representing the class $n.k.i$.  Under all possible labelings, there are $n!$ such polygons, many of which are identified in \M{n}.  Since the isotropy group maps the polygon to itself in the moduli space, we quotient $n!$ by the order of the isotropy group $|\IG_f|$, where $f$ is any face in \M{n} of the given class.  This provides another answer to for the Corollary above; combining the results, we find a solution to {\em P3}\, alternate to the one given in Proposition~\ref{p:cp3}.

\begin{prop}
The number of faces in $K_{n-1}$ of the class $n.k.i$ is 
$$\frac{n \cdot 2^{k+1}}{ |\IG_f|}.$$
\label{p:new}
\end{prop}


\subsection{}
The following sketches one possible algorithm for finding the group $\IG_f$.  This is done {\em locally} (looking at each cell of the cluster under twisting) and {\em globally} (looking at the entire dissected polygon under the dihedral group).

\begin{defn}
For a connected graph $G$, the {\em distance} between two vertices $x, y$ is the minimum of the lengths of paths containing $x$ and $y$. The {\em eccentricity} of a vertex is the maximum of its distances to all other vertices.  The {\em center} of $G$ is the subgraph induced by the vertices of minimum eccentricity.
\end{defn}

\noindent It is a classic result of graph theory that the center of a tree is either one vertex or one edge. Taking the dual, the center of the dissected polygon is either a cell or a diagonal.

\begin{figure}[h]
\centering {\includegraphics {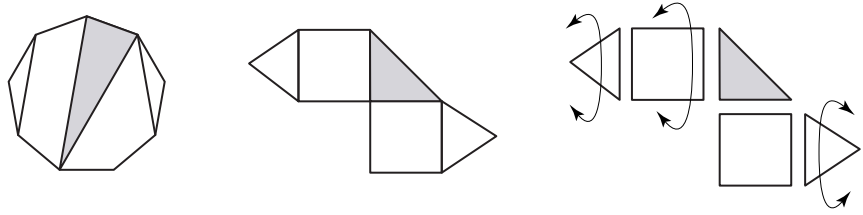}}
\caption{Local cells}
\label{localcell}
\end{figure}

\noindent {\em Non-central cells:} For each non-central cell of the polygon, observe if there exists an automorphism of the cell obtained by twisting along diagonals. If such a map exists, then this contributes $\Z_2$ to $\IG_f$.  Figure~\ref{localcell}b shows the class $9.4.10$ as a cluster, with the central cell shaded. Out of the four non-central cells, three of them have a $\Z_2$ action (Figure~\ref{localcell}c), contributing $\Z_2^3$ to the isotropy group.

\begin{figure}[h]
\centering {\includegraphics {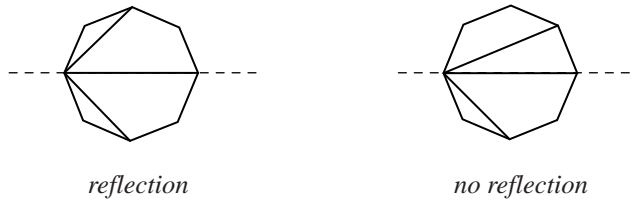}}
\caption{Central diagonal}
\label{centerdiag}
\end{figure}

\noindent {\em Central diagonal:} If the center is a diagonal, there are two possibilites: Observe whether or not there exists a (global) reflection along the central diagonal.  Figure~\ref{centerdiag} shows class $8.3.8$ with reflection and class $8.3.10$ without.

\begin{figure}[h]
\centering {\includegraphics {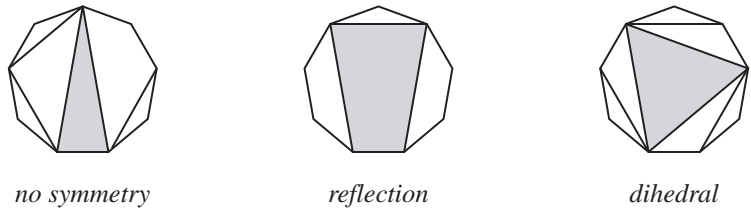}}
\caption{Central cell}
\label{centercell}
\end{figure}

\noindent {\em Central cell:} If the center is a cell, there are three possibilites.  Figure~\ref{centercell}a shows class $9.5.3$ with the central cell (shaded) having no symmetry. The cell could also have just a reflection action (class $9.3.11$), as in Figure~\ref{centercell}b. Finally, for an $n$-gon, there could be a dihedral action $D_k$ on the central cell, where $k$ divides $n$.  Figure~\ref{centercell}c shows class $9.6.4$ with the $D_3$ group action.\footnote{Although the picture shown has obvious $\Z_3$ rotation symmetry, reflection is also possible due to twisting.}

Since there are no intersecting diagonals, the different twists and reflections, along with any dihedral group action, all commute. Therefore, to find the isotropy group, we look at the  contribution of each local and global action and simply take their direct product.  For example, the isotropy groups (up to isomorphism) for the polygons of Figure~\ref{centercell} are $\Z_2^4$, $\Z_2^4$, and $\Z_2^3 \times D_3$, respectively.  Using Proposition~\ref{p:new}, we calculate $\trh$ for each of the classes to be $36$, $9$, and $24$, respectively, matching the values obtained from Proposition~\ref{p:cp3}.
 
\begin{figure}[h]
\centering {\includegraphics {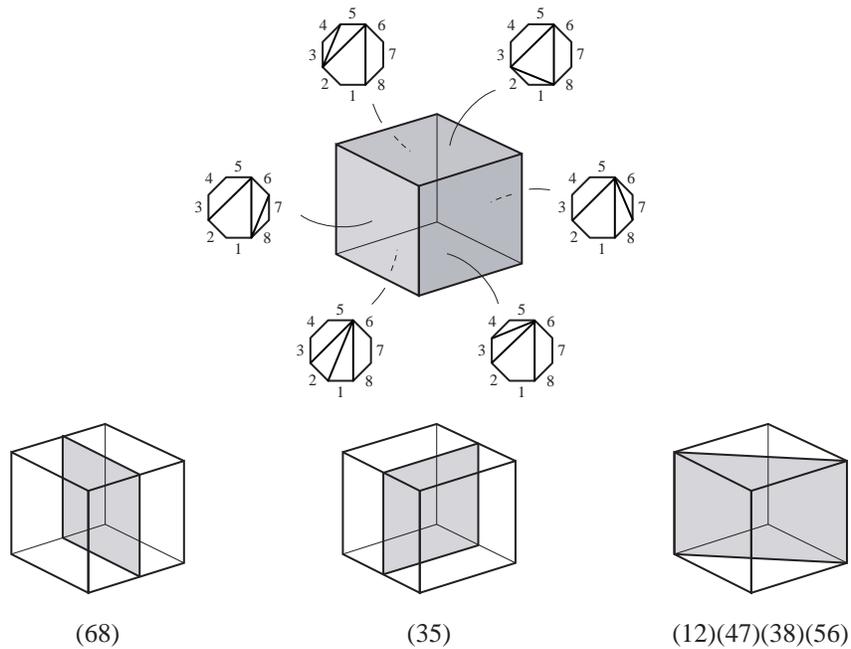}}
\caption{Labeled polytope}
\label{isocube}
\end{figure}

\begin{rem}
Each class corresponds to a polytope. By definition, each element of the isotropy group maps the polytope to itself; but what type of map is it? Using Figure~\ref{isotropy}a as an example of a face $f$, we show how the generators of $\IG_f$ act on the polytope.  Note that the polytope structure of $f$ is a cube, the product of three line segments, as labeled in Figure~\ref{isocube}. The actions of the three generators of $\IG_f$ given in Example~\ref{e:isotropy} are the reflections along the shaded planes.
\end{rem}

%
%

\section {Applications}

\subsection{}
\begin{figure}[h]
\centering {\includegraphics {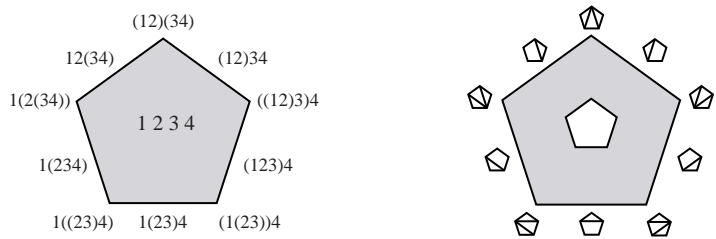}}
\caption{Parentheses and Polygons}
\label{k4}
\end{figure}

We first point out a few interesting similarities to classical combinatorics \cite[\S1.15]{comtet}.  As noted above, partitioning of an $n$-gon using $k$ diagonals corresponds to using $k$ sets of brackets meaningfully on $n-1$ {\em non-commutative} variables (see Figure~\ref{k4}). 
 
The {\em Generalized bracketing problem of Schr\"oder} is to find the number $s_n$ of different bracketings of $n$ variables using an arbitrary number of brackets (where each bracket must contain at least two elements).  Lemma~\ref{l:cayley} gives us the number of different bracketing for a fixed set of $k$ brackets.  Summing up, we find
$$s_n = \sum_{k=0}^{n-2} \frac{1}{k+1} \; \binom{n+k}{k} \; \binom{n-2}{k}.$$

\begin{figure}[h]
\centering {\includegraphics {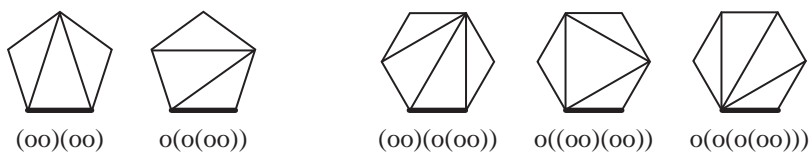}}
\caption{Representations of $w_4$ and $w_5$}
\label{f:wedder}
\end{figure}

The {\em Wedderburn-Etherington commutative bracketing problem} is to find the number $w_n$ of ways of using $n-2$ sets of brackets meaningfully on $n$ {\em commutative} variables. In other words, it is the number of different {\em bracketings} that arise from commuting binary subproducts. It is not too hard to see that $w_n$ gives the number of $(n+1)$-gons with maximal diagonals rooted at an outside edge, identified up to reflection along the root edge and twisting along diagonals, but {\em not} rotation.  It follows that $w_n = a_{n-1, n}$ (see Table~\ref{a-mn}).  Figure~\ref{f:wedder} shows $w_4 = 2$ and $w_5 = 3$ as polygon and bracketing representations.

\begin{rem}
The coefficients of the generating function of Equation~\ref{e:gena} generalize the Wedderburn-Etherington numbers:  The possible ways of using $k$ brackets on $n$ commutative variables is $a_{k+1, n}$.
\end{rem}


\subsection{}

We mention two occurrences outside mathematics that encode information using labeled trees.  The first is in biology; in particular, the areas of genetics and evolution \cite{ph}.  The theory of evolution conjectures that there exist links between certain species to common ancestors using the evidence from protein sequences.  An evolutionary tree, or {\em phylogeny}, is a means of organizing this data: The external (labeled) vertices represent species for which data is available whereas the internal (unlabeled) vertices are the hypothetical ancestors for which no direct evidence is known.  Phylogenies usually have no vertices of degree two.  The following emumerates the number of phylogenies having a fixed set of leaves.

\begin{prop} \textup{\cite[\S2]{fr}}
For $n$ external vertices, the number $T_n$ of labeled trees is given by the following recurrence relation, where $T_2=T_1=1$:
$$ T_n = (2-n) \, T_{n-1} \;+\; \sum_{1 \leq i \leq n-2} T_{i+1} \, T_{n-i} \, \binom{n-1}{i}.$$
\end{prop}

\begin{rem} A few observations follow:

1. Just as the generalized bracketing problem of Schr\"oder sums over all cells of $K_n$, the Proposition above sums over all the cells of \CM{n}.

2. A {\em rooted} phylogeny is a rooted labeled tree where the root does not receive a label. The root of a phylogeny represents the ancestor of all organisms of the tree. For $n$ labeled external vertices, the number $P_n$ of distinct rooted phylogenies is $P_n = T_{n+1}$ \cite[\S2]{fr}.

3. There exist many possible metrics on phylogenetic trees \cite{fr1}.  It was also noted in \cite{ph} that one can give lengths to the edges of labeled trees, making them {\em weighted}, which hold data on the time estimate of evolution taking place between two species of the tree.  Recently, a space of phylogenetic trees was constructed, closely related to the `dual' of \M{n} \cite{bhv}; the metric on this tree space inherits some properties from the natural metric on \M{n} \cite[\S2.2]{dev}.

4. The idea of phylogeny even extends to literature, in understanding the authenticity of the manuscripts of {\em The Canterbury Tales} \cite{ct}.
\qed
\end{rem}


\subsection{}
Another appearance of labeled trees is in physics; in particular, the quantum theory of angular momentum. Given a system of $n$ independent parts, each with its own angular momentum, the theory deals with constructing the total angular momentum by coupling the given $n$ momenta sequentially in pairs using $SU(2)$ Wigner coefficients. Biedenharn and Louck further relate this to the Racah-Wigner algebra and the 6-$j$ symbols \cite[\S5]{racah}. Enumerating the various {\em binary coupling schemes} gives a total of $c_n$ possibilities, thought of as counting maximal bracketings on $n$ variables (zero dim cells of $K_n$). Permuting two angular momenta implies a phase change in their {\em state vectors}.  Considering  binary coupling schemes over all such permutations  yield $n! \cdot c_n$ possibilites (zero dim cells in $n!$ copies of $K_n$).  However, in quantum physics it is not possible to distinguish state vectors differing by phase factors; schemes are identified by commuting binary subproducts. This yields $(2n-3)!!$ possible coupling schemes (zero dim cells in \M{n+1}).  As noted, the idea of binary coupling closely follows the gluing of associahedra to form the moduli space, at least in the zero dimensional case.  Although combinatorially the idea of coupling momenta can be generalized beyond binary pairs (by possibly counting the codim $k$ cells of moduli spaces), the meaning of this in the context of theoretical physics is not yet well understood.


\vspace{.5in}

{\small
\begin{ack}
The authors would like to acknowledge their indebtedness to, and admiration for, the work of Professor Arthur Cayley, who not only addressed a dissection problem similar to the one in this paper, but also laid the foundations for the enumerative study of tree-like structures in general.   This pioneering work is all the more remarkable in that it was only a small part of the mathematical output of this eminent Victorian.
\end{ack}}

%
%

\begin{rem}
For size limitations, we have not included Figures 22-26 in this paper.  One can obtain the paper {\bf with} the figures at:  {\em www.math.ohio-state.edu/$\sim$ devadoss}.
\end{rem}

\clearpage

%
%

\begin{table}[h]
{\scriptsize \begin{tabular}{|rrrrrrrrr|} \hline
3 & 4 & 5 & 6 & 7 & 8 & 9 & 10 & $Total$ \\ \hline
3 &&&&&&&& 5 \\ 
1 & 1 &&&&&&& 5 \\
&& 1 &&&&&& 1 \\ \hline
4 &&&&&&&& 14 \\ 
2 & 1 &&&&&&& 21 \\
1 && 1 &&&&&& 6 \\ 
& 2 &&&&&&& 3 \\
&&& 1 &&&&& 1 \\ \hline
5 &&&&&&&& 42 \\ 
3 & 1 &&&&&&& 84 \\
2 && 1 &&&&&& 28 \\
1 & 2 &&&&&&& 28 \\
& 1 & 2 &&&&&& 7 \\
1 &&& 1 &&&&& 7 \\
&&&& 1 &&&& 1 \\ \hline
6 &&&&&&&& 132 \\
4 & 1 &&&&&&& 330 \\
3 && 1 &&&&&& 120 \\
2 && 2 &&&&&& 180 \\
2 &&& 1 &&&&& 36 \\
1 & 1 & 1 &&&&&& 72 \\
& 3 &&&&&&& 12 \\
1 &&&& 1 &&&& 8 \\
& 1 && 1 &&&&& 8 \\
&& 2 &&&&&& 4 \\ 
&&&&& 1 &&& 1 \\ \hline
7 &&&&&&&& 429 \\
5 & 1 &&&&&&& 1287 \\ 
3 & 2 &&&&&&& 990 \\
3 &&& 1 &&&&& 165 \\
2 & 1 & 1 &&&&&& 495 \\
3 && 1 &&&&&& 495 \\
1 & 3 &&&&&&& 165 \\
2 &&&& 1 &&&& 45 \\
& 2 & 1 &&&&&& 45 \\
1 &&&&& 1 &&& 9 \\
& 1 &&& 1 &&&& 9 \\
&& 1 & 1 &&&&& 9 \\
&&&&&& 1 && 1 \\ \hline
8 &&&&&&&& 1430 \\ 
6 & 1 &&&&&&& 5005 \\
5 && 1 &&&&&& 2002 \\
4 & 2 &&&&&&& 5005 \\
4 &&& 1 &&&&& 715 \\
3 & 1 & 1 &&&&&& 2860 \\
2 & 3 &&&&&&& 1430 \\
3 &&&& 1 &&&& 220 \\
2 && 2 &&&&&& 330 \\
2 & 1 && 1 &&&&& 660 \\
1 & 2 & 1 &&&&&& 660 \\
& 4 &&&&&&& 55 \\ 
2 &&&&& 1 &&& 55 \\
1 & 1 &&& 1 &&&& 110 \\
1 && 1 & 1 &&&&& 110 \\
& 2 && 1 &&&&& 55 \\
& 1 & 2 &&&&&& 55 \\
1 &&&&&& 1 && 10 \\
& 1 &&&& 1 &&& 10 \\
&& 1 && 1 &&&& 10 \\
&&& 2 &&&&& 5 \\ 
&&&&&&& 1 & 1 \\ \hline
\end{tabular}}
\vspace{.2cm}
\caption{Dissections of the $n$-gon ($n=5,\ldots,10$)}
\label{t:diss}
\end{table}


\begin{table}[h]
{\scriptsize \begin{tabular}{|r|rrrrrrrrrrrr|} \hline
& 1 & 2 & 3 & 4 & 5 & 6 & 7 & 8 & 9 & 10 & 11 & 12 \\
2 & 1 &&&&&&&&&&& \\
3 & 1 & 1 &&&&&&&&&& \\
4 & 1 & 3 & 2 &&&&&&&&& \\
5 & 1 & 5 & 8 & 3 &&&&&&&& \\
6 & 1 & 8 & 22 & 20 & 6 &&&&&&& \\
7 & 1 & 11 & 46 & 73 & 49 & 11 &&&&&& \\
8 & 1 & 15 & 87 & 206 & 233 & 119 & 23 &&&&& \\
9 & 1 & 19 & 147 & 485 & 807 & 689 & 288 & 46 &&&& \\
10 & 1 & 24 & 236 & 1021 & 2320 & 2891 & 1988 & 696 & 98 &&& \\
11 & 1 & 29 & 356 & 1960 & 5795 & 9800 & 9737 & 5561 & 1681 & 207 && \\
12 & 1 & 35 & 520 & 3525 & 13088 & 28586 & 38216 & 31350 & 15322 & 4062 & 451 & \\
13 & 1 & 41 & 730 & 5989 & 27224 & 74280 & 127465 & 139901 & 97552 & 41558 & 9821 & 983 \\ \hline


\end{tabular}}
\vspace{.2cm}
\caption{$a_{m,n}$}
\label{a-mn}
\end{table}

\begin{table}[h]
{\scriptsize \begin{tabular}{|r|rrrrrrrrrrr|} \hline
& 1 & 2 & 3 & 4 & 5 & 6 & 7 & 8 & 9 & 10 & 11 \\
2 & 1 &&&&&&&&&& \\
3 & 1 & 2 &&&&&&&&& \\
4 & 1 & 5 & 5 &&&&&&&& \\
5 & 1 & 9 & 21 & 14 &&&&&&& \\
6 & 1 & 14 & 56 & 84 & 42 &&&&&& \\
7 & 1 & 20 & 120 & 300 & 330 & 132 &&&&& \\
8 & 1 & 27 & 225 & 825 & 1485 & 1287 & 429 &&&& \\
9 & 1 & 35 & 385 & 1925 & 5005 & 7007 & 5005 & 1430 &&& \\
10 & 1 & 44 & 616 & 4004 & 14014 & 28028 & 32032 & 19448 & 4862 && \\
11 & 1 & 54 & 936 & 7644 & 34398 & 91728 & 148512 & 143208 & 75582 & 16796 & \\
12 & 1 & 65 & 1365 & 13650 & 76440 & 259896 & 556920 & 755820 & 629850 & 293930 & 58786 \\ \hline


\end{tabular}}
\vspace{.2cm}
\caption{$v_{m,n}$}
\label{v-mn}
\end{table}

\begin{table}[h]
{\scriptsize \begin{tabular}{|r|rrrrrrrrrrrr|} \hline
& 1 & 2 & 3 & 4 & 5 & 6 & 7 & 8 & 9 & 10 & 11 & 12 \\
3 & 1 &&&&&&&&&&& \\
4 & 1 & 1 &&&&&&&&&& \\
5 & 1 & 2 & 2 &&&&&&&&& \\
6 & 1 & 3 & 7 & 4 &&&&&&&& \\
7 & 1 & 4 & 15 & 18 & 7 &&&&&&& \\
8 & 1 & 5 & 28 & 57 & 49 & 14 &&&&&& \\
9 & 1 & 6 & 45 & 138 & 196 & 123 & 29 &&&&& \\
10 & 1 & 7 & 69 & 288 & 601 & 626 & 313 & 60 &&&& \\
11 & 1 & 8 & 98 & 540 & 1533 & 2322 & 1899 & 778 & 127 &&& \\
12 & 1 & 9 & 136 & 943 & 3468 & 7095 & 8362 & 5565 & 1936 & 275 && \\
13 & 1 & 10 & 180 & 1544 & 7124 & 18813 & 29741 & 28350 & 15880 & 4776 & 598 & \\
14 & 1 & 11 & 235 & 2419 & 13635 & 44868 & 90869 & 115642 & 92210 & 44433 & 11777 & 1320 \\ \hline
\end{tabular}}
\vspace{.2cm}
\caption{$b_{m,n}$}
\label{b-mn}
\end{table}

\begin{table}[h]
{\scriptsize \begin{tabular}{|r|rrrrrrrrrrrrr|} \hline
& 1 & 2 & 3 & 4 & 5 & 6 & 7 & 8 & 9 & 10 & 11 & 12 & 13 \\
3 & 1 &&&&&&&&&&&& \\
4 & 1 & 1 &&&&&&&&&&& \\
5 & 1 & 1 & 1 &&&&&&&&&& \\
6 & 1 & 2 & 3 & 2 &&&&&&&&& \\
7 & 1 & 2 & 6 & 5 & 2 &&&&&&&& \\
8 & 1 & 3 & 11 & 17 & 12 & 4 &&&&&&& \\
9 & 1 & 3 & 17 & 37 & 44 & 23 & 6 &&&&&& \\
10 & 1 & 4 & 26 & 78 & 131 & 118 & 52 & 11 &&&&& \\
11 & 1 & 4 & 36 & 140 & 325 & 410 & 298 & 109 & 18 &&&& \\
12 & 1 & 5 & 50 & 248 & 728 & 1249 & 1279 & 766  & 244 & 37 &&& \\
13 & 1 & 5 & 65 & 396 & 1476 & 3246 & 4462 & 3763 & 1921 &  532 & 66 && \\
14 & 1 & 6 & 85 & 624 & 2811 &  7717 & 13497 & 15198 & 10920 & 4843 & 1196 & 135 & \\
15 & 1 & 6 & 106 & 929 & 5032 & 16773 & 36384 & 52041 & 49577 & 30848 & 12068 & 2671 & 265 \\
\hline
\end{tabular}}
\vspace{.2cm}
\caption{$f_{m,n}$}
\label{f-mn}
\end{table}

\clearpage

%
%

\bibliographystyle{amsplain}

\end{document}